\documentclass[a4paper]{article}
\begin{document}

\title{Geometric Algebra for Subspace Operations}
\author{T.A. Bouma, L. Dorst, H.G.J. Pijls\\{\it University of Amsterdam}}
\date{April 6th, 2001 CE}
\maketitle

\begin{abstract}
The set theory relations $\in$, $\backslash$, $\Delta$, $\cap$, and $\cup$ have
corollaries in subspace relations.  Geometric Algebra is introduced as the
ideal framework to explore these subspace operations.  The relations $\in$,
$\backslash$, and $\Delta$ are easily subsumed by Geometric Algebra for
Euclidean metrics.  A short computation shows that the meet ($\cap$) and join
($\cup$) are resolved in a projection operator representation with the aid of
one additional product beyond the standard Geometric Algebra products.  The
result is that the join can be computed even when the subspaces have a common
factor, and the meet can be computed without knowing the join.  All of the
operations can be defined and computed in any signature (including degenerate
signatures) by transforming the problem to an analogous problem in a different
algebra through a transformation induced by a linear invertible function (a LIFT
to a different algebra).  The new results, as well as the techniques by which
we reach them, add to the tools available for subspace computations.
\end{abstract}

\section{Introduction}
Operations on subspaces are useful in applications everywhere.  Geometric
Algebra, the intriguing algebra promoted by David Hestenes to unify and
simplify many areas of mathematics\cite{Hestenes-1,Hestenes-2,Hestenes-3}, is
introduced as the ideal framework to explore subspace operations.  A large
repertoire of operations to compute subspace operations are made available by
Geometric Algebra, but a few holes remain, notably with respect to the meet and
join of subspaces. This paper should resolve the outstanding issues.  A section
on preliminaries makes this treatment reasonably self-contained.  Four subspace
operations are introduced on equal footing, motivated by the set theory
relations ($\backslash$, $\Delta$, $\cap$, and $\cup$), common usage, and
applied needs.  An overview of the four subspace operations for all signatures
indicates that they are all fundamental and interconnected.  Experts will note
that Geometric Algebra deals with oriented subspaces and that in general there
is no oriented solution to the meet and join problem\cite{Stolfi-1}.  This
problem is avoided by representing unoriented subspaces by projection operators.
This allows us to extend the meet and join defined in the previous
literature\cite{Rota-1,Hestenes-2} to give meaningful (nonzero) results
for any subspaces.  The price we pay for this extension is that the meet and
join presented in this paper are not linear.  This is not a disadvantage,
because our meet and join agree with the previous literature except when the
previous literature results are zero.  Many operations in this paper have an
arbitrary scale and orientation.  This is addressed in the penultimate section,
where a geometrical significance can be given to the linear result of zero from
the previous literature.

\section{Preliminaries}
This section includes the definitions and motivation for four subspace
operations.  Following the new definitions is a review of Geometric Algebra.

\subsection{The Subspace Operations}
The set theory operations of $\backslash$, $\Delta$, $\cap$, and $\cup$ can be
applied to subspaces of $\mathcal{R}^n$.  However, in general the outcome will not
be a subspace, since the set theory operations do not respect the linear
structure of the subspaces.  Four subspaces operations are defined that are
motivated by the set theory operations, but that respect the linear structure
and thus always produce subspaces.  Let $\mathcal{A}$ and $\mathcal{B}$ be two
subspaces of $\mathcal{R}^n$.

\begin{enumerate}
\item
{\em The Meet Operation}\\
The meet of $\mathcal{A}$ and $\mathcal{B}$ is the set $\{\mathbf{x} \in
\mathcal{R}^n : \mathbf{x} \in \mathcal{A}$ and $\mathbf{x} \in \mathcal{B}\}$.
In words it is the largest common subspace.  It shall be denoted by
$\mathcal{A} \cap \mathcal{B}$.

\item
{\em The Join Operation}\\
The join of $\mathcal{A}$ and $\mathcal{B}$ is the set $\{\mathbf{x} \in
\mathcal{R}^n : \exists \mathbf{a} \in \mathcal{A}$ and $\exists \mathbf{b} \in
\mathcal{B}$ such that $\mathbf{x} = \mathbf{a} + \mathbf{b}\}$.  In words it
is the span of the two subspaces (i.e. the smallest common superspace).  It
shall be denoted by $\mathcal{A} \cup \mathcal{B}$.

\item
{\em The Difference Operation}\\
The difference of $\mathcal{A}$ and
$\mathcal{B}$ is the set $\{\mathbf{a} \in \mathcal{A}: \forall \mathbf{b} \in
\mathcal{B}\hspace{1em}\mathbf{a} \cdot \mathbf{b} = 0\}$, provided
$\mathcal{B}$ is a subspace of $\mathcal{A}$.  In words it is the orthogonal
complement of $\mathcal{B}$ in $\mathcal{A}$.  It shall be denoted by
$\mathcal{A} \backslash \mathcal{B}$.

\item
{\em The Symmetric Difference Operation}\\
The symmetric difference of
$\mathcal{A}$ and $\mathcal{B}$ is the set $\{\mathbf{x} \in \mathcal{R}^n:
\exists \mathbf{a} \in \mathcal{A}$ and $\exists \mathbf{b} \in \mathcal{B}$
such that $\mathbf{x} = \mathbf{a} + \mathbf{b}$ and $\forall \mathbf{c} \in
\mathcal{A} \cap \mathcal{B}\hspace{1em}\mathbf{x} \cdot \mathbf{c} = 0\}$.
In words it is the orthogonal complement of the meet in the join. It shall be
denoted by $\mathcal{A} \Delta \mathcal{B}$.  Clearly, 
$\mathcal{A} \Delta \mathcal{B} = (\mathcal{A} \cup \mathcal{B}) \backslash
(\mathcal{A} \cap \mathcal{B})$.
\end{enumerate}

This paper always uses the symbols $\backslash$, $\Delta$, $\cap$, and
$\cup$ to denote subspace operations.

\subsection{A Review of Geometric Algebra}

Geometric Algebra is based on Clifford Algebra over a real vector space.  A
Clifford Algebra is an algebra generated by the scalars and the elements of a
vector space with a metric.  The algebra has a linear, associative, distributive
product, and the square of a vector is the squared length determined
by the metric (for more details on Clifford Algebras see\cite{Lounesto-1}).  Let
$\mathcal{R}^{p,q,r}$ be a ($p$+$q$+$r$)-dimensional real vector space with a set of
linearly independent, mutually orthogonal vectors, $\{\mathbf{p}_1, ... ,
\mathbf{p}_p,\mathbf{q}_1, ... ,\mathbf{q}_q, \mathbf{r}_1, ... , \mathbf{r}_r\}$,
such that $\mathbf{p}_i^2=1$, $\mathbf{q}_j^2=-1$, and $\mathbf{r}_k^2=0$.  The
Clifford Algebra and the Geometric Algebra generated by $\mathcal{R}^{p,q,r}$ shall be
denoted $\mathcal{R}_{p,q,r}$. 

The essential difference between a Geometric Algebra and a Clifford Algebra is
that the elements of a Geometric Algebra are given a geometric interpretation.
This leads to a focus on operations that are defined on geometrically meaningful
subsets of the algebra and therefore to the introduction of additional structure
(and more products between elements), so that a consistent geometric
interpretation
can be maintained on the results of computations.  To emphasize this difference,
we call the elements of the Geometric Algebra {\em multivectors} and we call the
standard (Clifford) product between elements in the Algebra the {\em geometric
product}.  The geometric product is denoted by juxtaposition of operands, as in
$AB$.  

A summary of the extra features and terminology of Geometric Algebra follows.

\begin{enumerate}
\item
{\em Blades}\\
If a nonzero multivector can be written as the geometric product of $r$
mutually  anticommuting vectors, then it is called an $r$-blade.  The word blade
refers to an $r$-blade with the value of $r$ unspecified.  Real numbers are
considered $0$-blades and are often called scalars.  Vectors are considered
$1$-blades.  The square of a blade is a scalar.  For a vector the square is the
squared length determined by the metric.  Zero is considered an $r$-blade for
any value of $r$, but any nonzero blade is an $r$-blade for only one value of
$r$\cite{Hestenes-1}.  

\item
{\em Steps (or Grades)}\\
A linear combination of $r$-blades will be called an $r$-vector, and will be
said to have {\em step} (or grade) $r$.  The space of $r$-vectors is a linear
subspace of the entire Clifford Algebra.  An arbitrary multivector, $A$, can be
uniquely written as $\sum \langle A \rangle_r$ where $\langle A \rangle_r$ is
the $r$-vector part of $A$ if $r$ is a nonnegative integer and $\langle A
\rangle_r$ is zero if $r$ is not a nonnegative integer.

\item
{\em Outer Product}\\
The outer product of an $r$-vector $A$ and an $s$-vector $B$ is defined to be
$\langle AB \rangle_{s+r}$.  It is denoted by $A \wedge B$ and is extended by
linearity to arbitrary multivectors.  The outer product is associative between
all multivectors and anti-symmetric between vectors.  The outer product of two
blades is a blade (see the appendix).

\item
{\em (Contraction) Inner Product}\\
The (contraction) inner product of an $r$-vector $A$ and an $s$-vector $B$ is
defined to be $\langle AB \rangle_{s-r}$. It is denoted by $A \rfloor B$ and is
extended by linearity to arbitrary multivectors.  This inner product differs
slightly from the inner product of Hestenes\cite{Hestenes-1}.  It has the useful
properties that $(A \wedge B) \rfloor C = A \rfloor (B \rfloor C)$ and
$\mathbf{a}A = \mathbf{a} \rfloor A + \mathbf{a} \wedge A$ for any multivectors
$A$, $B$, $C$ and any vector $\mathbf{a}$. These details are obvious from
considering this definition and looking at the proofs in
Hestenes\cite{Hestenes-1}.  The inner product is explicitly expanded for the
vectors $\mathbf{a}$, $\mathbf{c}_1$, $\mathbf{c}_2$, ... , and $\mathbf{c}_k$ as

\begin{equation}\label{inner-expansion}
\mathbf{a} \rfloor (\mathbf{c}_1 \wedge \mathbf{c}_2 \wedge ... \wedge
\mathbf{c}_k) = \sum (-1)^{r+1} (\mathbf{a} \rfloor \mathbf{c}_i)
\mathbf{c}_1 \wedge ... \wedge \check{\mathbf{c}}_r \wedge ... \wedge
\mathbf{c}_k 
\end{equation}

where the inverted circumflex indicates that the $r$th vector was omitted from
the product.  The proof in reference\cite{Hestenes-1} carries over to the
contraction inner product with no modification.  The inner product of two
blades is a blade (see the appendix).

\item
{\em Subspaces}\\
One of the geometric interpretations common in geometric algebra is to use
blades to represent subspaces.  This works because blades are closely related
to subspaces.  Given a nonzero blade, $\mathbf{A}$, the set $\{\mathbf{x} \in
\mathcal{R}^n : \mathbf{x} \wedge \mathbf{A} =0\}$ is a subspace of $\mathcal{R}^n$.
Similarly given an oriented basis \{$\mathbf{x}_1, \mathbf{x}_2, ... , 
\mathbf{x}_k$\} for a $k$-dimensional subspace $\mathcal{A}$, there is a
nonzero blade that corresponds to that oriented basis.  If $k=0$ then that
blade is $1$.  If $k \neq 0$ then that blade is $\mathbf{x}_1 \wedge
\mathbf{x}_2 \wedge ... \wedge \mathbf{x}_k$.  The identity $\mathbf{x} \in
\mathcal{A} \iff \mathbf{x} \wedge \mathbf{A} = 0$ is the means by which
Geometric Algebra subsumes the operation $\in$ from set theory.  Since the
square of a blade is a scalar, the inverse of a blade (if it has one), is equal
to the blade divided by its square, so the inverse is just a scalar multiple of
the original blade, and hence a blade and its inverse represent the same
(unoriented) subspace.

\item
{\em Pseudoscalars}\\
Given an algebra, $\mathcal{R}_{p,q,r}$, and a set of $(p+q+r)$ vectors
$\{\mathbf{v}_1, ... , \mathbf{v}_{p+q+r}\}$, the outer product
$\mathbf{v}_1 \wedge ... \wedge \mathbf{v}_{p+q+r}$ is called
a pseudoscalar.  Often a particular nonzero pseudoscalar is
singled out.  This preferred pseudoscalar serves to determine
the reference orientation for the vector space and sometimes is
used to perform duality operations.  Usually the preferred
pseudoscalar is called {\em the} pseudoscalar and denoted
$\mathbf{I}$.  This can be seen as merely a special case of the
previous section on subspaces, by noticing that every nonzero
$k$-blade is a pseudoscalar for the $k$ dimensional subspace it
represents.

\item
{\em Projection Operators}\\
For a Euclidean metric, blades are closely related to projection operators.
Given a nonzero blade, $\mathbf{A}$, representing the subspace $\mathcal{A}$,
the vector $P_{\mathbf{A}}(\mathbf{x})=(\mathbf{x} \rfloor \mathbf{A})
\mathbf{A}^{-1}$ is the orthogonal projection of the vector $\mathbf{x}$ onto
the subspace $\mathcal{A}$. As taken from\cite{Hestenes-1}, the identities

\begin{equation}\label{projection-plus}
\mathbf{AB} = \mathbf{A} \wedge \mathbf{B} \Rightarrow P_{\mathbf{AB}} =
P_{\mathbf{B}} + P_{\mathbf{A}}
\end{equation}

and

\begin{equation}\label{projection-minus}
\mathbf{AB} = \mathbf{A} \rfloor \mathbf{B} \Rightarrow P_{\mathbf{AB}} =
P_{\mathbf{B}} - P_{\mathbf{A}}
\end{equation}

hold for any nonzero blades $\mathbf{A}$ and $\mathbf{B}$.  Equation
(\ref{projection-plus}) implies that a projection operator can be decomposed
analogously to the way its corresponding blades can be factored.  

\item
{\em Outermorphism}\\
Given a linear function, $f : \mathcal{R}^{p,q,r} \rightarrow \mathcal{R}^{p,q,r}$, there is
an extension of the function to arbitrary multivectors called the outermorphism
of $f$. It is denoted by $\underline{f}$ and is defined to be the identity when
restricted to the scalars.  The condition $\underline{f}$($A \wedge (B + C)$) =
$\underline{f}$($A$)$ \wedge \underline{f}$($B$) + $\underline{f}$($A$)$ \wedge
\underline{f}$($C$) is then sufficient to define the outermorphism on arbitrary
multivectors.  This extension is well-established\cite{Hestenes-3}.

\end{enumerate}

A notation convention is adopted to aid the reader in easily making
meaningful distinctions between different multivectors.
Lowercase Greek letters are reserved for scalars.  Lowercase Latin letters
are reserved for integers or functions when not in bold face, while
lowercase Latin letters are reserved for vectors when in bold face.  Bold
face is reserved for blades.  Lastly, uppercase Latin letters are
used when it is impossible or unnecessary to be more specific about the nature
of a multivector.  This notation convention simplifies the reading of equations
and emphasizes that different geometric interpretations are applied to
different elements of the algebra.

\section{The Euclidean Metrics}
Since every blade represents an oriented subspace and blades are easy to
compute with, they are a natural candidate for subspace computations.
The extra scalar degree of freedom allows the future potential for more precise
calculations with subspaces that attach meaning to the magnitude of a blade.
Therefore in this paper four blade operations are introduced to correspond to
the four subspace operations.  The blade operations are made first for Euclidean
space because the relationship between blades and projection operators is
strongest in Euclidean space.  Therefore first a correspondence is made from
projection operators to blades, then the four blade operations are defined and
shown to faithfully mirror the subspace operations.

\subsection{The Blade Correspondence}

Here we give an algorithm to construct a blade from its corresponding projection
operator.  The algorithm has an arbitrary scale and orientation inherited from
an arbitrarily chosen basis, which is the best that can be expected.  Let
$P_{\mathcal{A}}$ be an idempotent linear operator on $\mathcal{R}^n$ such that the
image of $\mathcal{R}^n$ is a $t$-dimensional subspace, $\mathcal{A}$.  The algorithm
constructs a blade that characterizes $\mathcal{A}$ as follows.  First we
construct a set of candidate blades, and then show that all the nonzero
candidate blades represent the subspace $\mathcal{A}$.  Finally we show
that at least one of the candidate blades is, in fact, nonzero.  Let
$\mathbf{V}_1$, $\mathbf{V}_2$, ... , $\mathbf{V}_k$ be $k$ $t$-blades such
that $\{\mathbf{V}_1$, $\mathbf{V}_2$, ... ,$\mathbf{V}_k$\} is a
basis for the space of $t$-vectors.  Let $\mathbf{T}_i = 
\underline{P_{\mathcal{A}}}(\mathbf{V}_i)$. The set $\{\mathbf{T}_1,
\mathbf{T}_2, ... , \mathbf{T}_k\}$ is our set of candidate blades.  By the
properties of the outermorphism, each $\mathbf{T}_i$ is a $t$-blade.  Each
$\mathbf{T}_i$ is clearly the outer product of $t$ vectors and each of these $t$
vectors is in $\mathcal{A}$.  If the $t$ vectors are linearly
dependent then $\mathbf{T}_i = 0$, if not then they form a basis for
$\mathcal{A}$ so $\mathbf{T}_i \neq 0$ and $\mathbf{T}_i$ is exactly the kind
of blade to characterize the subspace $\mathcal{A}$.  All that remains is to
show that at least one of the candidate blades is nonzero.  Since a blade
that characterizes $\mathcal{A}$ exists we know that it is a linear combination
of \{$\mathbf{V}_1$ , $\mathbf{V}_2$, ... , $\mathbf{V}_k$\}, so by the
linearity of the outermorphism there must be a $\mathbf{V}_i$ such that
$\mathbf{T}_i \neq 0$.

Given this correspondence we can translate operations between projection
operators into operations between blades, except for a loss of the scale and
orientation.

\subsection{The Inner Division Operation}

Consider two nonzero blades, $\mathbf{A}$ and $\mathbf{B}$, that characterize
the subspaces $\mathcal{A}$ and $\mathcal{B}$ respectively.  When $\mathcal{B}$
is a subspace of $\mathcal{A}$ we use the expression $\mathbf{A} \backslash
\mathbf{B}$ to denote the quantity $\mathbf{B}^{-1} \mathbf{A}$ and the call
the operation {\em inner division}.  The inner division operation is
motivated by the difference operation for subspaces.  The justification
requires showing two points, first that under such conditions, $\mathbf{B}^{-1}
\mathbf{A}$ is a blade, and second that $\mathbf{x} \wedge (\mathbf{B}^{-1}
\mathbf{A}) = 0 \iff \mathbf{x} \in \mathcal{A} \backslash \mathcal{B}$.

Let $\{\mathbf{b}_1, ... , \mathbf{b}_s\}$ be an orthogonal basis for
$\mathcal{B}$.  Let $\{\mathbf{a}_1, ... , \mathbf{a}_r\}$ be an orthogonal
basis for $\mathcal{A} \backslash \mathcal{B}$.  Clearly $\{\mathbf{b}_1, ... ,
\mathbf{b}_s, \mathbf{a}_1, ... , \mathbf{a}_r\}$ is an orthogonal basis for
$\mathcal{A}$.  Since $\{\mathbf{b}_1, ... , \mathbf{b}_s\}$ is an orthogonal
basis for $\mathcal{B}$ and since $\mathbf{B}^{-1}$ characterizes the subspace
$\mathcal{B}$ it follows that $\mathbf{B}^{-1}$ is a nonzero scalar multiple of
$\mathbf{b}_1\mathbf{b}_2 ... \mathbf{b}_s$.  Similarly $\mathbf{A}$ is a
nonzero scalar multiple of $\mathbf{b}_1\mathbf{b}_2 ... \mathbf{b}_s\mathbf{a}_1
\mathbf{a}_2 ... \mathbf{a}_r$.  It then follows that there exists two nonzero
scalars $\alpha$ and $\beta$ such that $\mathbf{B}^{-1} \mathbf{A} = \alpha
(\mathbf{b}_1\mathbf{b}_2 ... \mathbf{b}_s)(\mathbf{b}_1\mathbf{b}_2 ...
\mathbf{b}_s\mathbf{a}_1\mathbf{a}_2 ... \mathbf{a}_r) = \beta \mathbf{a}_1
\mathbf{a}_2 ... \mathbf{a}_r$.  Thus $\mathbf{B}^{-1} \mathbf{A}$ is a blade
and it characterizes the subspace $\mathcal{A} \backslash \mathcal{B}$.

A quick look at the step of the output reveals that when $\mathcal{B}$ is a
subspace of $\mathcal{A}$ then $\mathbf{B}^{-1} \mathbf{A} =
\mathbf{B}^{-1} \rfloor \mathbf{A}$.  Therefore, while it appears that the
inner division operation is based on the geometric product it is also just
as easily based on the inner product.  It is useful whenever a product between
two blades can be written as either of two products, because either definition
can be used from line to line of a computation, depending on which product gives
simplifications at that particular moment.  An example is the identity,

\begin{equation}\label{inner-decomposition}
\mathbf{A} = \mathbf{B}(\mathbf{A} \backslash \mathbf{B}) = \mathbf{B} \wedge
(\mathbf{A} \backslash \mathbf{B})
\end{equation}
which is proved by decomposing the inner division first into the geometric
product and then into the inner product and looking at the step of the outcome.

\subsection{The Delta Product}
Consider two nonzero blades, $\mathbf{A}$ and $\mathbf{B}$, that characterize
the subspaces $\mathcal{A}$ and $\mathcal{B}$ respectively.  Let $\mathcal{C} =
\mathcal{A} \cap \mathcal{B}$ and let $\mathbf{C}$ be any blade characterizing
that subspace.  When $\mathcal{C} = \{0\}$, $\mathbf{A} \wedge \mathbf{B} \neq
0$. When $\mathcal{C} \neq \{0\}$, $\mathbf{A} \wedge \mathbf{B} = 0$.  In the
latter case we can define $\mathbf{A}_{\perp}=\mathbf{AC}$ and
$\mathbf{B}_{\perp}=\mathbf{C}^{-1}\mathbf{B}$.  Since $\mathbf{C}$ and
$\mathbf{C}^{-1}$ both represent $\mathcal{C}$, which is a subspace of both
$\mathcal{A}$ and $\mathcal{B}$, the previous section makes it clear that
$\mathbf{A}_{\perp}$ and $\mathbf{B}_{\perp}$ are blades.  The intersection
of the subspaces characterized by $\mathbf{A}_{\perp}$ and $\mathbf{B}_{\perp}$
contains only the element zero, so $\mathbf{A}_{\perp} \wedge
\mathbf{B}_{\perp} \neq 0$.

By construction $\mathbf{AB} = \mathbf{A}_{\perp}\mathbf{B}_{\perp}$, so the
highest step portion of $\mathbf{AB}$ is $\mathbf{A}_{\perp} \wedge
\mathbf{B}_{\perp}$, and therefore a blade.  This motivates a new product for
blades which we call the {\em delta product}.  The delta product of two blades,
$\mathbf{A}$ and $\mathbf{B}$ is denoted $\mathbf{A} \Delta \mathbf{B}$ and
defined to be the highest step portion of $\mathbf{AB}$.  The delta product is
motivated by the symmetric difference operation for subspaces.  The
justification requires showing that $\mathbf{x} \wedge (\mathbf{A} \Delta
\mathbf{B}) = 0 \iff \mathbf{x} \in \mathcal{A} \Delta \mathcal{B}$.

Let $\{\mathbf{a}_1, ... , \mathbf{a}_r\}$ be an orthogonal basis for
$\mathcal{A} \backslash \mathcal{C}$.  Let $\{\mathbf{b}_1, ... ,
\mathbf{b}_s\}$ be an orthogonal basis for $\mathcal{B} \backslash \mathcal{C}$.
Similar arguments as in the previous section demonstrate that there exists a
nonzero scalar $\alpha$ such that $\mathbf{AB} = \alpha\mathbf{a}_1\mathbf{a}_2
... \mathbf{a}_r\mathbf{b}_1\mathbf{b}_2 ... \mathbf{b}_s$, so $\mathbf{A} \Delta
\mathbf{B} = \alpha \mathbf{a}_1 \wedge \mathbf{a}_2 \wedge ... \wedge
\mathbf{a}_r \wedge \mathbf{b}_1 \wedge \mathbf{b}_2 \wedge ... \wedge
\mathbf{b}_s$.  Clearly $\mathbf{x} \wedge (\mathbf{A} \Delta \mathbf{B}) = 0
\Rightarrow \mathbf{x} \in \mathcal{A} \Delta \mathcal{B}$. Therefore assume
$\mathbf{x} \in \mathcal{A} \Delta \mathcal{B}$ and we will show that
$\mathbf{x} \wedge (\mathbf{A} \Delta \mathbf{B}) = 0$.  By definition $\exists
\mathbf{a} \in \mathcal{A}$ and $\exists \mathbf{b} \in \mathcal{B}$ such that
$\mathbf{x} = \mathbf{a} + \mathbf{b}$ and $\forall \mathbf{c} \in \mathcal{C}
\hspace{1em}\mathbf{x} \cdot \mathbf{c} = 0$.  Let $\mathbf{a}_0 = \mathbf{a} +
P_{\mathcal{C}}(\mathbf{b})$ and $\mathbf{b}_0 = \mathbf{b} + P_{\mathcal{C}}(
\mathbf{a})$.  Clearly $\mathbf{a}_0 + \mathbf{b}_0 = \mathbf{x} + P_{
\mathcal{C}}(\mathbf{x})=\mathbf{x}$.  Clearly $\mathbf{a}_0 \in \mathcal{A}$
and $\mathbf{a}_0 \rfloor \mathbf{C} = \mathbf{a} \rfloor \mathbf{C} -
\mathbf{b} \rfloor \mathbf{C}=0$, so $\forall \mathbf{c} \in \mathcal{C}$,
$\mathbf{c} \rfloor \mathbf{a}_0 = 0$ therefore $\mathbf{a}_0 \in \mathcal{A}
\backslash \mathcal{C}$.  This means $\mathbf{a}_0 \wedge \mathbf{a}_1 \wedge
\mathbf{a}_2 \wedge ... \wedge \mathbf{a}_r = 0$ so $\mathbf{a}_0 \wedge
(\mathbf{A} \Delta \mathbf{B})=0$.  Similarly for $\mathbf{b}_0$.  Therefore
$\mathbf{x} \wedge (\mathbf{A} \Delta \mathbf{B})=0$.

The delta product is very different than the inner or outer product in its
algebraic properties.  The major difference is that $\mathbf{A} \Delta
(\mathbf{B} + \mathbf{C}) \neq \mathbf{A} \Delta \mathbf{B} + \mathbf{A} \Delta
\mathbf{C}$, so the delta product cannot be extended by linearity to arbitrary
multivectors.  The delta product can only be used on blades.  More care must
be taken with implementations of the delta product because even a small change
in either $\mathbf{A}$ or $\mathbf{B}$ can cause a change in the step of
$\mathbf{A}\Delta\mathbf{B}$.

\subsection{The Meet Operation}
Consider two nonzero blades, $\mathbf{A}$ and $\mathbf{B}$, that characterize
the subspaces $\mathcal{A}$ and $\mathcal{B}$ respectively.  The blade
correspondence can be used to define a new product for blades called the {\em
meet} and denoted $\mathbf{A} \cap \mathbf{B}$.  $\mathbf{A} \cap \mathbf{B}$ is
defined modulo a scale and an orientation as the blade corresponding to the
projection operator $P_{\mathbf{A} \cap \mathbf{B}}$, where $P_{\mathbf{A} \cap
\mathbf{B}}$ is defined as:

\begin{equation}\label{Meet}
P_{\mathbf{A} \cap \mathbf{B}}(\mathbf{x})=\frac{P_{\mathbf{B}}
(\mathbf{x})-P_{\mathbf{A} \Delta \mathbf{B}}(\mathbf{x})+P_{(\mathbf{A} \Delta
\mathbf{B}) \mathbf{B}^{-1}}(\mathbf{x})}{2}
\end{equation}

The justification requires that $P_{\mathbf{A} \cap \mathbf{B}}$ be the
orthogonal projection onto $\mathcal{A} \cap \mathcal{B}$.  Let $\mathbf{C}$ be
a blade characterizing the subspace $\mathcal{A} \cap \mathcal{B}$, then define
$\mathbf{A}_{\perp}=\mathbf{AC}$ and $\mathbf{B}_{\perp}=\mathbf{C}^{-1}
\mathbf{B}$ as above.  Now the result follows from a simple application of
equation (\ref{projection-plus}).  First note that $\mathbf{B}_{\perp}=
\mathbf{B}\backslash\mathbf{C}$, so equation (\ref{inner-decomposition}) implies
that $\mathbf{B}=\mathbf{C}\mathbf{B}_{\perp}=\mathbf{C} \wedge \mathbf{B}_{
\perp}$, therefore we have the following identity: 

\begin{equation}\label{B}
P_{\mathbf{B}} = P_{\mathbf{B}_{\perp}}+P_{\mathbf{C}}
\end{equation}

Similarly $\mathbf{A}_{\perp} \wedge
\mathbf{B}_{\perp}=\mathbf{B}_{\perp}((\mathbf{A}_{\perp} \wedge
\mathbf{B}_{\perp})\backslash \mathbf{B}_{\perp})=\mathbf{B}_{\perp} \wedge
((\mathbf{A}_{\perp} \wedge \mathbf{B}_{\perp})\backslash \mathbf{B}_{\perp})$,
therefore we have the following identity:

\begin{equation}\label{A Delta B}
P_{\mathbf{A} \Delta \mathbf{B}} = P_{\mathbf{B}_{\perp}^{-1}(\mathbf{A} \Delta
\mathbf{B})} + P_{\mathbf{B}_{\perp}}
\end{equation}

Finally $(\mathbf{A} \Delta \mathbf{B})\mathbf{C}^{-1}=(\mathbf{A} \Delta
\mathbf{B}) \wedge \mathbf{C}^{-1}$, because every vector in $\mathbf{C}^{-1}$
is orthogonal to every vector in $\mathbf{A} \Delta \mathbf{B}$.  Since
$(\mathbf{A} \Delta \mathbf{B})\mathbf{B}_{\perp}^{-1}$ represents a subspace of
$\mathbf{A} \Delta \mathbf{B}$ it is just as true that $((\mathbf{A} \Delta
\mathbf{B})\mathbf{B}_{\perp}^{-1})\mathbf{C}^{-1}=((\mathbf{A} \Delta
\mathbf{B})\mathbf{B}_{\perp}^{-1}) \wedge \mathbf{C}^{-1}$, therefore we have
the following identity:

\begin{equation}\label{(A Delta B)B}
P_{(\mathbf{A} \Delta \mathbf{B})\mathbf{B}^{-1}} = P_{(\mathbf{A} \Delta
\mathbf{B})\mathbf{B}_{\perp}^{-1}} + P_{\mathbf{C}^{-1}}
\end{equation}

The linear combination of the three projection operators has now been {\em
reduced} to the linear combination of four projection operators.  A quick appeal
to equation (\ref{projection-minus}) implies that for two blades, if their
geometric product is a scalar then their projection operators are equal.  Now
$\mathbf{C}\mathbf{C}^{-1}$ is a scalar, and so is $\mathbf{B}_{\perp}^{-1}
(\mathbf{A} \Delta \mathbf{B})(\mathbf{A} \Delta \mathbf{B})
\mathbf{B}_{\perp}^{-1}$.  Therefore the terms $P_{(\mathbf{A} \Delta
\mathbf{B})\mathbf{B}_{\perp}^{-1}}$ and $P_{\mathbf{B}_{\perp}^{-1}(\mathbf{A}
\Delta \mathbf{B})}$ cancel and the terms $P_{\mathbf{C}^{-1}}$ and
$P_{\mathbf{C}}$ combine.  The result, (equation (\ref{Meet})), then follows
from equations (\ref{B}), (\ref{A Delta B}), and (\ref{(A Delta B)B}).

A small commentary is in order.  The first comment is that the calculation of
the blade correspondence will be simplified by the fact that if the steps of
$\mathbf{A}$, $\mathbf{B}$, and $\mathbf{A} \Delta \mathbf{B}$ are $r$, $s$, 
and $q$ respectively, then the step of the meet is $\frac{r+s-q}{2}$.  The
second comment is that the blade correspondence is not precise about the scale
and orientation of the blade because the blade correspondence inherits an
arbitrary scale and orientation from an arbitrary basis of blades.   Since the
meet for blades is defined by the blade correspondence, this lack of precision is
then passed on to the meet for blades, except for the {\em disjoint} case.  The
disjoint case occurs when $\mathcal{A} \cap \mathcal{B} = \{ 0 \}$, and in this
case one can choose a basis {\em a priori}.  This is possible because in this
case the meet for blades is a scalar.  Therefore one can choose the scalar `$1$'
for the basis, and then since the blade correspondence uses the outermorphism of
the projection operator and an outermorphism is the identity when restricted to
the scalars, the blade correspondence gives a determinate answer for the meet,
namely it gives `$1$' back again.  The third comment is that the meet for blades
given here is different from previous literature\cite{Rota-1,Hestenes-2}, which
only relates nontrivially to our definition when $\mathcal{A} \cup \mathcal{B} =
\mathcal{R}^n$.  When $\mathcal{A} \cup \mathcal{B} \neq \mathcal{R}^n$ the previous
literature gives the zero blade as the result, while we can also treat that
case.  The price we pay is linearity.  Like the delta product, the meet is not
linear. 

\subsection{The Join Operation}
Previously we noted that $(\mathbf{A} \Delta \mathbf{B})\mathbf{C}=(\mathbf{A}
\Delta \mathbf{B}) \wedge \mathbf{C}$, in fact $(\mathbf{A} \Delta \mathbf{B})
\wedge \mathbf{C}$ characterizes $\mathcal{A} \cup \mathcal{B}$, therefore using
equation (\ref{projection-plus}), we find that $P_{\mathbf{A} \cup \mathbf{B}}=
P_{\mathbf{A} \Delta \mathbf{B}}+P_{\mathbf{A} \cap \mathbf{B}}$.  We have two
alternatives.  The first alternative is to define the projection operator for
the join directly as

\begin{equation}
P_{\mathbf{A} \cup \mathbf{B}}(\mathbf{x})=\frac{P_{(\mathbf{A} \Delta 
\mathbf{B})\mathbf{B}^{-1}}(\mathbf{x})+P_{\mathbf{A} \Delta 
\mathbf{B}}(\mathbf{x})+P_{\mathbf{B}}(\mathbf{x})}{2}
\end{equation}
and if the steps of $\mathbf{A}$, $\mathbf{B}$, and $\mathbf{A} \Delta
\mathbf{B}$ are $r$, $s$, and $q$ respectively, then the step of the join is
$\frac{r+s+q}{2}$.  The second alternative is to define the join for blades
directly in terms of the meet for blades and the inner division operation
through the equation $\mathbf{A} \cup \mathbf{B} = \mathbf{A} \wedge
(\mathbf{B} \backslash (\mathbf{A} \cap \mathbf{B}))$.  Just as the meet had a
definite scale and orientation only in the disjoint case, this definition
allows the join to inherit the definite scale and orientation from the meet,
since in that case $\mathbf{A} \cup \mathbf{B} = \mathbf{A} \wedge \mathbf{B}$.
It bears mentioning that this join for blades only agrees with the previous
literature\cite{Rota-1,Hestenes-2} in the disjoint case, but again
definitions in the previous literature are merely zero when $\mathcal{A} \cap
\mathcal{B} \neq \{0\}$, so this definition is an extension.  As with the
meet the price we pay is linearity.

\section{The Non-Euclidean Metrics}
A new concept is introduced to allow the subspace operations to be performed
in any metric.  Each of the four operations is then investigated in turn.

\subsection{A LIFT Between Clifford Algebras}
Given two algebras, $\mathcal{R}_{p,q,r}$ and $\mathcal{R}_{a,b,c}$, such that $p$+$q$+$r$=
$a$+$b$+$c$=$n$, and a linear invertible function, $f$, from the vectors of
$\mathcal{R}_{p,q,r}$ to $\mathcal{R}_{a,b,c}$ then $\underline{f}$ is a linear invertible
map between the two algebras.  Call the extended function $\underline{f}$ a LIFT
(`linear invertible function' transformation) of $\mathcal{R}_{p,q,r}$ to
$\mathcal{R}_{a,b,c}$.

A LIFT can be used to transfer a problem with subspaces to another algebra,
preserving incidence relations but allowing the metric to change.  Often a LIFT 
is taken to a Euclidean space.  After solving the problem in that space,
subspace results can be pulled back to the original space.  Examples that extend
the previous results on the subspace operations follow.

A minor variation of the LIFT is for $f$ to go into an $n$-dimensional subspace
of a Clifford Algebra over a larger vector space, then the outermorphism is a
linear invertible map onto a subalgebra.  This is especially nice when the
outermorphism is an isomorphism between the original algebra and the
subalgebra.  Such a LIFT is called an embedding LIFT (or e-LIFT).  This is
especially common for degenerate algebras $\mathcal{R}_{p,q,r}$, for which an e-LIFT
to $\mathcal{R}_{p+r,q+r,0}$ always exists.  To see the existence of the e-LIFT, let
$\{\mathbf{p}_1, ... ,\mathbf{p}_p,\mathbf{q}_1, ... ,\mathbf{q}_q,\mathbf{r}_1,
... , \mathbf{r}_r\}$ be an orthogonal basis for the vectors in $\mathcal{R}_{p,q,r}$
such that $\mathbf{p}_i^2=1$, $\mathbf{q}_j^2=-1$, and $\mathbf{r}_k^2=0$ and
let $\{\mathbf{e}_1, ... ,\mathbf{e}_{p+r},\mathbf{f}_1, ... ,\mathbf{f}_{q+r}\}$
be an orthogonal basis for the vectors in $\mathcal{R}_{p+r,q+r,0}$ such that
$\mathbf{e}_i^2=1$ and $\mathbf{f}_j^2=-1$.  Then let $f$ be a linear
function such that $f(\mathbf{p}_i)=\mathbf{e}_i$, $f(\mathbf{q}_j)=
\mathbf{f}_j$, and $f(\mathbf{r}_k)=\mathbf{e}_{p+k}+\mathbf{f}_{q+k}$.
Then $\underline{f}$ is the promised isomorphism.

\subsection{The Meet Operation}
If one fixes an arbitrary LIFT, $\underline{f}$, from $\mathcal{R}_{p,q,r}$ to
$\mathcal{R}_{p+q+r,0,0}$ then the meet, $\mathbf{A} \cap \mathbf{B}$, between two
blades $\mathbf{A}$ and $\mathbf{B}$ can be defined by:

\begin{equation}
\mathbf{A} \cap \mathbf{B} = \underline{f}^{-1}(\underline{f}(\mathbf{A}) \cap
\underline{f}(\mathbf{B}))
\end{equation}

The scale and orientation of $\underline{f}(\mathbf{A}) \cap \underline{f}(
\mathbf{B})$ are indeterminate except when $\underline{f}(\mathbf{A})$ and
$\underline{f}(\mathbf{B})$ are disjoint, which only happens when $\mathbf{A}$
and $\mathbf{B}$ are disjoint.  The LIFT is an outermorphism, so it is the
identity on the scalars, so the meet has a definite scale and orientation in the
disjoint case, and they are the same scale and orientation as in the Euclidean
space.  More importantly, in the disjoint case, the scale and orientation are
independent of which LIFT, $\underline{f}$, was chosen.  The preservation of the
outer product and the scalars makes it clear that this meet corresponds to
$\mathcal{A} \cap \mathcal{B}$.  This means that this definition has a
well-defined scale and orientation in exactly the cases where the Euclidean
definition did, and it has an arbitrary scale and orientation in exactly the
cases where the Euclidean definition did.

\subsection{The Join Operation}
If one fixes fixes an arbitrary LIFT, $\underline{f}$, from $\mathcal{R}_{p,q,r}$ to
$\mathcal{R}_{p+q+r,0,0}$ then the join, $\mathbf{A} \cup \mathbf{B}$, between two
blades $\mathbf{A}$ and $\mathbf{B}$ can be defined by:

\begin{equation}
\mathbf{A} \cup \mathbf{B} = \underline{f}^{-1}(\underline{f}(\mathbf{A}) \cup
\underline{f}(\mathbf{B}))
\end{equation}

The scale and orientation of $\underline{f}(\mathbf{A}) \cup \underline{f}(
\mathbf{B})$ are indeterminate except when $\underline{f}(\mathbf{A})$ and
$\underline{f}(\mathbf{B})$ are disjoint, which only happens when $\mathbf{A}$
and $\mathbf{B}$ are disjoint, in which case the join should reduce to the outer
product.  The LIFT is an outermorphism, so it preserves the outer product, so
clearly $\underline{f}^{-1}(\underline{f}(\mathbf{A}) \wedge \underline{f}(
\mathbf{B}))=\underline{f}^{-1}(\underline{f}(\mathbf{A} \wedge \mathbf{B}))=
\mathbf{A} \wedge \mathbf{B}$.  Therefore the join has a definite scale and
orientation in the disjoint case, and the scale and orientation are independent
of which LIFT, $\underline{f}$, was chosen.  The preservation of the outer
product and the scalars makes it clear that this join corresponds to
$\mathcal{A} \cup \mathcal{B}$.  This means that this definition has a
well-defined scale and orientation in exactly the cases where the Euclidean
definition did, and it has an arbitrary scale and orientation in exactly the
cases where the Euclidean definition did.

\subsection{The Inner Division Operation}
Consider two nonzero blades, $\mathbf{A}$ and $\mathbf{B}$ in $\mathcal{R}_{p,q,r}$,
that characterize the subspaces $\mathcal{A}$ and $\mathcal{B}$ respectively
such that $\mathcal{B} \subseteq \mathcal{A}$.  When $\mathbf{B}^{-1}$ exists
we can calculate $\mathbf{A} \backslash \mathbf{B} =\mathbf{B}^{-1} \mathbf{A}$
as usual. Otherwise, we need an e-LIFT, $f$, from $\mathcal{R}_{p,q,r}$ to
$\mathcal{R}_{p+r,q+r,0}$.  Let $\mathbf{I}$ be the pseudoscalar of
$\mathcal{R}_{p+r,q+r,0}$.  Then define $\mathbf{A} \backslash \mathbf{B}=
\underline{f}^{-1}(\underline{f}(\mathbf{A}) \cap (\underline{f}(
\mathbf{B})I))$.  This meets the definition for the subspace operation,
but now the scale and orientation has an arbitrary dependence on $f$.

\subsection{The Delta Product}
Consider two nonzero blades, $\mathbf{A}$ and $\mathbf{B}$ in $\mathcal{R}_{p,q,r}$,
that characterize the subspaces $\mathcal{A}$ and $\mathcal{B}$ respectively.
If $\mathbf{A} \cap \mathbf{B}$ has an inverse then the geometric product,
$\mathbf{AB}$, is nonzero and the highest step portion represents the symmetric
difference, $\mathbf{A} \Delta \mathbf{B}$ as usual.  Otherwise, we need an
e-LIFT, $f$, from $\mathcal{R}_{p,q,r}$ to $\mathcal{R}_{p+r,q+r,0}$.  Let $\mathbf{I}$ be
the pseudoscalar of $\mathcal{R}_{p+r,q+r,0}$.  Then define $\mathbf{A} \Delta
\mathbf{B} = \underline{f}^{-1}((\underline{f}(\mathbf{A}) \cup \underline{f}(
\mathbf{B})) \cap ((\underline{f}(\mathbf{A}) \cap \underline{f}(
\mathbf{B}))I))$.  This meets the definition for the subspace operation, but
now the scale and orientation has an arbitrary dependence on $f$.  Note that
the symmetric difference was used to compute the meet and join for Euclidean
signatures, but the meet and join are used to define the symmetric difference
for Non-Euclidean signatures.

\section{Linearity and the Meaning of Zero}
The meet and join for blades presented in this paper are not linear, (e.g. in
general $(\mathbf{A}+\mathbf{B}) \cap \mathbf{C} \neq \mathbf{A} \cap \mathbf{C}
+\mathbf{B} \cap \mathbf{C}$, even when $\mathbf{A}+\mathbf{B}$ is a nonzero
blade).  The previous literature\cite{Rota-1,Hestenes-2} have linear results for
the meet and join.  In our notation the meet and join of the previous literature
are $(A\mathbf{I}) \rfloor B$ and $A \wedge B$ respectively.

The linear versions can operate on any multivector (by linear extension), but
the geometric interpretation of the computation becomes confusing.  Also, even
when the linear versions operate on nonzero blades, they can disagree with the
subspace operations by giving a result of zero.  It is not surprising that the
subspace operations disagree with the linear operations, because it was exactly
the preponderance of the answer $0$ for many meaningful computations that
motivated the creation of the new subspace operations of this paper.  However
the opposing versions can be reconciled by pursuing a geometric interpretation
of the zero blade.

An interpretation of the zero blade that is consistent with the general
enterprise of representing oriented subspaces by blades is to have the zero
blade represent an {\em indeterminate} oriented subspace. At first glance, it
appears that the zero blade represents the whole space (i.e. $\{\mathbf{x} \in
\mathcal{R}^n : \mathbf{x} \wedge 0=0\}=\mathcal{R}^n$), but this interpretation would imply
that the zero blade represents a different subspace depending on a particular
enveloping pseudoscalar (a property that destroys natural subalgebra and
enveloping algebra relationships).  But since this is implicit in stating that
the zero blade can represent {\em any} subspace this is actually support for the
interpretation proposed here.  Furthermore, since the hope of representing
subspaces by blades is to eventually be able to deal with uncertainty in
geometrical computations, the scale factor would naturally be used to represent
how well-determined the blade is.  This implies that the zero blade represents a
completely undetermined subspace. Lastly the linear versions give the zero blade
as the result if either of the input blades is the zero blade, this adsorbing
property is consistent with the indeterminacy the zero blade represents.

When two oriented subspaces do not span the pseudoscalar, the linear meet gives
the result of the zero blade because the orientation of the meet cannot be
determined.  The linear meet is a fully functional quantitative operation, which
can give a quantitative meet, but only {\em if} given a quantitative join first
(in the role of the pseudoscalar).  Similarly when two oriented subspaces have a
nontrivial intersection, the linear join gives the result of the zero blade
because the join has an undetermined orientation.

With this interpretation for the zero blade, the linear meet and join can be
compared to the versions presented in this paper.  The linearity can be an 
advantage for implementation for some applications, and if that advantage
outweighs the costs of getting the zero blade as a result, then an educated
decision to implement the linear versions can be made for that application.
The interpretation for the zero blade can also be used to extend the subspace
operations for nonzero blades to the zero blade by declaring the inner division,
the delta product, the meet, and the join to be zero if either input of the two
input blades is the zero blade.

\section{Conclusion}
Blades can represent subspaces, and in applications we need to perform
operations on subspaces.  Therefore we would naturally want operations on
blades that mirror the results of computations that we would have liked to
perform on subspaces.  The four subspace operations (inner division, delta
product, meet, and join) are the first four operations from nonzero blades to
nonzero blades.  The hope is that these operations can contribute to
quantitative computations with oriented subspaces.

The four operations are on equal footing because the delta product and the
inner division are used to define the meet and join in Euclidean signatures,
while the meet and join are used to define the delta product and inner division
in Non-Euclidean signatures.  This indicates that they are all fundamental
(though not as fundamental as the geometric product) and that they are all
interconnected.

Standard concepts in geometric algebra needed to be augmented because the meet
and join for blades cannot be defined\cite{Stolfi-1} with an orientation due to
fundamental geometric problems.  This fundamental problem was solved by using
projection operators to represent unoriented subspaces.  Within this solution
the delta product helps to compute the meet and join for blades.

The price to be paid for this augmentation is that the new blade operations are
not linear and cannot be extended to arbitrary multivectors.  The authors
believe that their non-linearity might be the reason that the operations have
not been used previously.  It is only by sacrificing linearity, and thus losing
applicability to arbitrary multivectors, that one can solve the meet and join
for blades.

The four subspace operations are tools intended for general use in applications,
however no applied examples are included in this paper.  Readers looking for
examples of the inner division and the delta product need look no farther than
the proof of the meet for Euclidean metrics, and readers looking for examples of
the meet and join need look no farther than the inner division and delta product
for Non-Euclidean metrics.  Beyond the four subspace operations, this paper
utilizes another tool of general applicability. This tool is the LIFT
(`linear invertible function' transformation).  It is an invertible map between
algebras of different signatures.  This tool can be used to advantageously
transform problems that are independent of signature to whichever signature is
most helpful at any particular moment.  This tool is demonstrated in the paper
by extending the results of the meet and join from the Euclidean case to the
Non-Euclidean case (even to degenerate signatures).

Lastly some loose ends are resolved.  A geometric interpretation is given to the
zero blade that explains the results of the previous literature and facilitates
the educated choice between different versions of the subspace operations.
The final loose end is resolved by the appendix, which includes proofs to
demonstrate that the inner and outer products go from blades to blades (even in
degenerate signatures).

\section{Appendix}

In this Appendix we prove that the outer product of two blades
is a blade and that the inner product of two blades is a blade.  To our
surprise, this does not appear to have been proved before, but is less trivial
than may have been assumed when degenerate algebras are considered.

\subsection{The Outer Product}

By associativity of the outer product it suffices to show that the outer
product of a vector and an $r$-blade is a blade.  The result is trivial because
the ($r$+1)-vector determines an ($r$+1)-dimensional subspace and that subspace
has an orthogonal basis.  However a more constructive proof is desirable, to
assist in the proof for the inner product and to see how such a factorization
can be made. The proof proceeds by induction on $r$, the base case is
$r=1$, or the outer product of two vectors is a blade.

Let $\mathbf{v}_1$ and $\mathbf{v}_2$ be two vectors.  Since $\frac{1}{2}
(\mathbf{v}_1 - \mathbf{v}_2)(\mathbf{v}_1 + \mathbf{v}_2) = \mathbf{v}_1
\wedge \mathbf{v}_2 + \frac{\mathbf{v}^2_1 - \mathbf{v}^2_2}{2}$, this gives a
factorization of $\mathbf{v}_1 \wedge \mathbf{v}_2$ when $\mathbf{v}^2_1 =
\mathbf{v}^2_2$.  If $\mathbf{v}^2_1 \neq \mathbf{v}^2_2$ then either
$\mathbf{v}_1^2$ or $\mathbf{v}_2^2$ is not equal to $0$.  Since $\mathbf{v}_1
\wedge \mathbf{v}_2 = -\mathbf{v}_2 \wedge \mathbf{v}_1$, we may assume
without loss of generality that $\mathbf{v}_1^2 \neq 0$.  Then we note that
$\mathbf{v}_1^{-1} \rfloor (\mathbf{v}_1 \wedge \mathbf{v}_2)$ is a vector and
that $\mathbf{v}_1(\mathbf{v}_1^{-1} \rfloor (\mathbf{v}_1 \wedge
\mathbf{v}_2)) = \mathbf{v}_1 \wedge \mathbf{v}_2$, so this gives a
factorization of $\mathbf{v}_1 \wedge \mathbf{v}_2$.  Therefore in both cases
the outer product of two vectors can be factored.

Assume that the outer product of $r$ vectors is a blade.  Let $\{\mathbf{a}, 
\mathbf{c}_1, ... ,\mathbf{c}_r\}$ be $r$+$1$ vectors.  The inductive step has
three cases.
\begin{enumerate}
\item
The first case is when $\mathbf{a} \rfloor \mathbf{c}_i = 0$ for
each $i$ .  Then by equation (\ref{inner-expansion}), $\mathbf{a} \wedge
\mathbf{c}_1 \wedge ... \wedge \mathbf{c}_r = \mathbf{a} (\mathbf{c}_1 \wedge
... \wedge \mathbf{c}_r)$.  By the inductive hypothesis, $\mathbf{c}_1 \wedge
... \wedge \mathbf{c}_r$ is a blade, so there exist  $r$ anticommuting vectors,
$\{\mathbf{a}_1, ... , \mathbf{a}_r\}$, such that $\mathbf{a}_1\mathbf{a}_2 ...
\mathbf{a}_r=\mathbf{c}_1 \wedge ... \wedge \mathbf{c}_r$.  Each $\mathbf{a}_i$
is a linear combination of the set $\{\mathbf{c}_1, ... ,\mathbf{c}_r\}$, so each
$\mathbf{a}_i$ anticommutes with $\mathbf{a}$, so $\mathbf{a}\mathbf{a}_1 ...
\mathbf{a}_r$ is a factorization of $\mathbf{a} \wedge \mathbf{c}_1 \wedge ...
\wedge \mathbf{c}_r$, therefore $\mathbf{a} \wedge \mathbf{c}_1 \wedge ...
\wedge \mathbf{c}_r$ is an ($r$+$1$)-blade.
\item
The next case is when $\mathbf{a}^2 \neq 0$.  Let $\mathbf{b}_i =
\mathbf{a}^{-1} (\mathbf{a} \wedge \mathbf{c}_i)$.  Then $\mathbf{a} \wedge
\mathbf{c}_1 \wedge ... \wedge \mathbf{c}_r = \mathbf{a} \wedge \mathbf{b}_1
\wedge ... \wedge \mathbf{b}_r$.  Now we have guaranteed that
$\mathbf{a} \rfloor \mathbf{b}_i = 0$ for each $i$, so the previous case
resolves the factorization.
\item
The last case is when there exists a $k$ such that
$\mathbf{a} \rfloor \mathbf{c}_k \neq 0$ and $\mathbf{a}^2 = 0$.  Without loss
of generality we may assume $k=1$ since the order of the vectors only affects
the sign of the outcome.  By the base case, $\mathbf{a} \wedge \mathbf{c}_1$ is
a $2$-blade.  Direct computation shows that the square of $\mathbf{a} \wedge
\mathbf{c}_1$ is $(\mathbf{a} \rfloor \mathbf{c}_1)^2$, hence nonzero, therefore
there exist two invertible vectors $\mathbf{c}$ and $\mathbf{d}_1$ such that
$\mathbf{a} \wedge \mathbf{c}_1 = \mathbf{c}\mathbf{d}_1$.  Let $\mathbf{d}_i =
(\mathbf{c}\mathbf{d}_1)^{-1}((\mathbf{c}\mathbf{d}_1) \wedge \mathbf{c}_{i})$
for $i \in \{2, ... , r\}$.  Then $\mathbf{a} \wedge \mathbf{c}_1 \wedge ...
\wedge \mathbf{c}_r = \mathbf{c} \wedge \mathbf{d}_1 \wedge \mathbf{d}_2 \wedge
... \wedge \mathbf{d}_r$.  Now since $\mathbf{c}^2 \neq 0$ the previous case
resolves the factorization.
\end{enumerate}

\subsection{The (Contraction) Inner Product}

Let $\mathbf{A} = \mathbf{a}_1 \wedge \mathbf{a}_2 \wedge ... \wedge
\mathbf{a}_m$ and $\mathbf{B}$ be two blades in $\mathcal{R}_{p,q,r}$.  By the
properties of the inner product, $\mathbf{A} \rfloor \mathbf{B} = \mathbf{a}_1
\rfloor ( \mathbf{a}_2 \rfloor ( ... \rfloor (\mathbf{a}_m \rfloor \mathbf{B})
... ))$, so it suffices to show that $\mathbf{a} \rfloor \mathbf{C}$ is a
blade for any vector $\mathbf{a}$ and any blade $\mathbf{C}$.  If $\mathbf{C}$
is a scalar then $\mathbf{a} \rfloor \mathbf{C} = 0$ so it is a blade.
Otherwise let $\mathbf{C} = \mathbf{c}_1 \wedge \mathbf{c}_2 \wedge ... \wedge
\mathbf{c}_k$.  The nonzero blade $\mathbf{C}$ characterizes a subspace with
signature $(s,t,u)$ and hence a subalgebra $\mathcal{R}_{s,t,u}$.  Using equation
(\ref{inner-expansion}), it is clear that $\mathbf{a} \rfloor \mathbf{C}$
resides in $\mathcal{R}_{s,t,u}$.  Let $\underline{f}$ be a LIFT  from $\mathcal{R}_{s,t,u}$
to $\mathcal{R}_{k,0,0}$. Since $\underline{f}(X)$ is a blade if and only if $X$ is a
blade then to show that $\mathbf{a} \rfloor \mathbf{C}$ is a blade it suffices
to show that $\underline{f}(\mathbf{a} \rfloor \mathbf{C})$ is a blade.  From
equation (\ref{inner-expansion}), it follows that $\underline{f}(\mathbf{a}
\rfloor \mathbf{C})$ is the sum of ($k$-$1$)-blades.  Let $\mathbf{I}$ be a
nonzero pseudoscalar for $\mathcal{R}_{k,0,0}$.  Let $\{\mathbf{H}_1, ... ,
\mathbf{H}_m\}$ be $m$ ($k$-$1$)-blades in $\mathcal{R}_{k,0,0}$.  Then $\mathbf{H}_1
+ ... + \mathbf{H}_m = (\mathbf{H}_1 + ... + \mathbf{H}_m)\mathbf{I}^{-1}
\mathbf{I} = \mathbf{v}\mathbf{I}$, where $\mathbf{v} = \mathbf{H}_1
\mathbf{I}^{-1} + ... + \mathbf{H}_m \mathbf{I}^{-1} = \mathbf{H}_1 \rfloor
\mathbf{I}^{-1} + ... + \mathbf{H}_m \rfloor \mathbf{I}^{-1}$ is a vector. If
$\mathbf{v} = 0$ then the sum, $\mathbf{H}_1 + ... + \mathbf{H}_m = \mathbf{v}
\mathbf{I}$, is zero and trivially a blade.  If not then $\mathbf{v} \neq 0$
and therefore $\mathbf{v}$ has an inverse and by the details of the proof for
the outer product, it is clear that $\mathbf{I}$ can be factored under the
geometric product with $\mathbf{v}^{-1}$ as a factor.  Then $\mathbf{H}_1 +
... + \mathbf{H}_m = \mathbf{v}\mathbf{I}$ is clearly a blade.

\end{document}